\documentclass{amsart}

\usepackage{amsmath,amsthm,amssymb,amscd,enumerate,mathtools,enumitem}
\usepackage{amsfonts}
\usepackage{rotating}
\usepackage{euscript}
\usepackage{pst-node}
\usepackage{epsfig,verbatim}
\usepackage{stackengine,scalerel}
\stackMath

\def\be{\begin{equation}}
\def\ee{\end{equation}}

\def\C{{\mathbb C}} 
\def\f{\EuScript}
 
\def\P{{\mathbb P}}

\def\Q{{\mathbb Q}}

\def\phi{{\varphi}}
\def\v{{\varepsilon}} 
\def\tt{\widetilde}
\def\deg{{\rm deg\,}}

\def\GCD{{\rm GCD }}

\def\bp{\begin{proposition}}
\def\ep{\end{proposition}}

\def\bt{\begin{theorem}}
\def\et{\end{theorem}}
\def\br{\begin{remark}}
\def\er{\end{remark}}
\def\be{\begin{equation}}
\def\bee{\begin{equation*}}
\def\l{\label}

\def\ee{\end{equation}}
\def\eee{\end{equation*}}
\def\bl{\begin{lemma}}
\def\el{\end{lemma}}
\def\bc{\begin{corollary}}
\def\ec{\end{corollary}}
\def\pr{\noindent{\it Proof. }}

\def\bd{\begin{definition}}
\def\ed{\end{definition}}
\def\t{\widetilde}

\def\hat{\widehat}
\def\Aut{{\rm Aut}}

\newtheorem{theorem}{Theorem}[section]
\newtheorem{lemma}[theorem]{Lemma}
\newtheorem{definition}[theorem]{Definition}
\newtheorem{corollary}[theorem]{Corollary}
\newtheorem{proposition}[theorem]{Proposition}
\newtheorem{problem}[theorem]{Problem}

\theoremstyle{definition}

\theoremstyle{definition}
\newtheorem{remark}[theorem]{Remark}
\def\bpr{\begin{problem}}
\def\epr{\end{problem}}

\mathtoolsset{showonlyrefs}

\begin{document}

\title{
Algebraic functions with  infinitely many values 
in a number field
}

\author[F. Pakovich]{Fedor Pakovich}
\thanks{
This research was supported by ISF Grant  No. 1092/22}
\address{Department of Mathematics, Ben Gurion University of the Negev, Israel}
\email{
pakovich@math.bgu.ac.il}

\begin{abstract}
We describe algebraic curves \( X : F(x, y) = 0 \) defined over \(\overline{\mathbb{Q}}\) that satisfy the following property: there exist a number field \(k\) and an infinite set \(S \subset k\) such that, for every \(y \in S\), the roots of the polynomial \(F(x, y)\) belong to \(k\).

\end{abstract}

%\keywords{}
%\subjclass[2010]{}

\maketitle

\section{Introduction}
The Hilbert irreducibility theorem, in its simplest form, states that 
if \( F(x, y) \) is an irreducible polynomial over \( \mathbb{Q} \), 
then for infinitely many integers \( n \), the polynomial \( F(x, n) \) 
remains irreducible in \( \mathbb{Q}[x] \). 
This result has numerous generalizations, the most relevant for this paper 
being the following extension to number fields: if \( F(x, y) \) is 
irreducible over a number field \( k \), then for infinitely many elements 
\( \tau \) in the ring of integers \( \mathcal{O}_k \), the polynomial 
\( F(x, \tau) \) is irreducible in \( k[x] \).

The Hilbert irreducibility theorem is equivalent to the statement that, 
for infinitely many values of \( n \), the degree of the field 
\( \mathbb{Q}(x_n) \), generated by a root \( x_n \) of \( F(x, n) \), 
attains the maximum possible degree over \( \mathbb{Q} \). 
However, this theorem provides no information about the fields 
\( \mathbb{Q}(x_n) \) themselves. In particular, the following question 
is outside the scope of the Hilbert theorem: how many distinct fields can be 
among \( \mathbb{Q}(x_n) \)? A closely related question is: how quickly do 
the degrees of fields generated by adjoining a sequence of such roots 
\( x_n \) to \( \mathbb{Q} \) grow?

This latter problem was investigated by Dvornicich and Zannier in \cite{z1}, 
where they proved that if \( \deg_x F \geq 2 \), then there exist constants 
\( c \) and \( N_0 \), depending on \( F \), such that 
\( \mathbb{Q}(x_1, x_2, \dots, x_N) \) has degree at least 
\( e^{\frac{cN}{\log N}} \) over \( \mathbb{Q} \) whenever \( N \geq N_0 \). 
Dvornicich and Zannier's result was later generalized to number fields 
by Bilu (\cite{bilu}), who showed that if \( F(x, y) \) is defined over a 
number field \( k \) of degree \( d \), then the field obtained by adjoining 
to \( k \) roots \( x_{\tau} \) of \( F(x, \tau) \), where \( \tau \) 
ranges over all elements of \( \mathcal{O}_k \) with height less than \( B \), 
has degree at least \( e^{\frac{cB^d}{\log B}} \) for \( B \) greater than \( B_0 \).
For further results in this direction, see \cite{bl1}, \cite{bl2}, and \cite{z2}.

If, instead of the sequence of fields  
$ \mathbb{Q}(x_1, x_2, \dots, x_N) $, one considers the sequence  
$ \mathbb{Q}(x_{i_1}, x_{i_2}, \dots, x_{i_N}) $, where  
$i_1, i_2, \ldots, i_N, \ldots$ is an arbitrary infinite sequence of integers, the degrees of  
$ \mathbb{Q}(x_{i_1}, x_{i_2}, \dots, x_{i_N}) $ over $\mathbb{Q}$ do not necessarily grow.  
For example, if $F(x, y) = x^r - y$, then for any $y$ that is an $r$th power of an integer,  
the roots of $F(x, y)$ lie in the same field $k$ generated by the $r$th roots of unity.  

This phenomenon motivates the following general question:  
for which irreducible algebraic curves  
$X : F(x, y) = 0$ over $\overline{\mathbb{Q}}$ do there exist a number field $k$  
and an infinite set $S \subset k$ such that, for each $y \in S$, the $y$-fiber of $X$,  
defined as the set of roots of $F(x, y)$, is contained in $k$?

A related question was studied by Zannier \cite{za}, who imposed the stronger condition that
$S$ lies in the ring of integers $\mathcal{O}_k$.
If, for such a set $S$, the corresponding $y$-fibers of a curve $X$ lie in a number field $k$,
then Siegel’s theorem implies that the Galois closure $\overline{\Sigma}$ of the field extension
$\Sigma = \mathbb{C}(X)/\mathbb{C}(y)$ has genus $0$, and 
Zannier’s results describe the possible fields $k$ and the Galois groups
$\mathrm{Gal}(\overline{\Sigma}/\Sigma)$.
We also mention the work of Bilu \cite{bilu0}, whose results imply that the {\it finite} sets  
$S \subset \mathcal{O}_k$ satisfying the above property can be described effectively,  
which is essentially equivalent to an effective version of Siegel’s theorem for integral points  
on Galois coverings of the projective line (see also the Appendix to \cite{z1}).

For brevity, we will refer to curves for which there exists a number field $k$
and an infinite set $S \subset k$ such that, for each $y \in S$, the $y$-fiber of $X$ is contained in $k$
as {\it curves with infinitely many $y$-fibers in a number field}.
In this paper, we describe such curves and, in particular, show
that the problem is equivalent to a seemingly stronger version—
in the spirit of arithmetic dynamics—asking whether there exists a set $S$ as above
that forms the forward orbit of a point under a rational function. 
Notice that 
any graph \( x = P(y) \), where 
\( P \in \overline{\mathbb{Q}}(z) \), obviously  has infinitely 
many \( y \)-fibers in a number field. Thus, we will always 
assume that the curves we consider satisfy \( \deg_x F \geq 2 \).

Our main result is the following statement.

\bt \l{t1}
  Let $X: F(x,y)=0$ be an  irreducible affine curve over $\overline\Q$ such that $\deg_xF\geq 2.$ Then the following conditions are equivalent. 

\begin{enumerate} [label=\arabic*)]

\item The curve $X$ has infinitely many $y$-fibers in a number field.  
\item 
There exists a rational function $A$ of degree at least two  over $\overline\Q$ satisfying the following condition:  for every  $y_0\in \overline\Q$,  
there exists a number field $k$ such that $y$-fibers of $X$ belong to $k$ for every \( y \in \overline\Q\) that belongs to the forward orbit of $y_0$ under $A$.

\item  
There exists a self-rational map \(\psi : X \to X\) of degree at least two, 
defined over \(\overline{\mathbb{Q}}\), which maps each \(y\)-fiber of \(X\) 
bijectively onto another \(y\)-fiber for all but finitely many values 
\(y \in \mathbb{C}\).

\item The Galois closure of the field extension $\C(X)/\C(y)$ has genus zero or one. 

\end{enumerate}
\et 
Notice that the first condition of the theorem readily follows from the second.  
It is also easy to verify that a self-rational map satisfying the third condition can be written in the form  
\be \l{po}  
\psi : (x, y) \mapsto (\phi(x, y), A(y)),  
\ee  
where $A$ is a rational function of one variable satisfying the second condition.  
Hence, the third condition entails the second.  
Finally, the theorem of Faltings yields that the first condition leads to the fourth.  
Thus, the essential step in the proof of Theorem \ref{t1} is to establish the direction $4) \Rightarrow 3)$,  
whose proof relies on a description of semiconjugacies between holomorphic maps  
of compact Riemann surfaces obtained in a series of papers \cite{semi}, \cite{arn}, \cite{rec}, \cite{lattes}.

Let us mention that the fourth condition of Theorem \ref{t1} 
can be effectively verified solely in terms of the ramification of \( y \), 
and that the function \( A \) in the second condition  and 
the  map \( \psi \) in the third condition (more precisely, families of such functions and maps) can also be effectively constructed.

We illustrate Theorem \ref{t1} with the following two examples (for more details,  
see Section \ref{ex}). Let $X : F(x, y) = 0$ be an algebraic curve of genus zero.  
Then $X$ admits a generically one-to-one parametrization $z \mapsto (U(z), V(z))$  
by some rational functions $U$ and $V$, and the fourth condition of the theorem is satisfied, for example, for $V = 3z^4 - 4z^3$ and arbitrary $U$.  
A possible choice of $A$ in this case is the rational function  
\be \l{aa} 
A = \frac{z^2 \left( z^3 - 240z^2 + 19200z - 512000 \right)}{625z^4 + 16000z^3 + 153600z^2 + 655360z + 1048576}.
\ee
Indeed, one can check that for the function  
\be \l{bb} 
B = -\frac{z^2 \left( 3z^3 - 10z^2 + 20z - 40 \right)}{15z^4 - 20z^3 + 32},
\ee
the diagram  
\be \l{dd} 
\begin{CD}
\C\P^1 @> B >> \C\P^1 \\
@V {V} VV @VV {V} V \\ 
\C\P^1 @> A >> \C\P^1
\end{CD}
\ee 
commutes, and $B$ sends the roots of $V(z) = z_0$ to the roots of $V(z) = A(z_0)$  
generically bijectively, since $\C(V,B) = \C(z)$.  
In turn, this implies that the diagram  
\be 
\begin{CD}
\C\P^1 @> B >> \C\P^1 \\
@V {(U,V)} VV @VV {(U,V)} V \\ 
X @> \psi >> X
\end{CD}
\ee 
commutes for some  $\psi$ of the form \eqref{po}, which  
sends $y$-fibers of $X$ to $y$-fibers.   
In particular, if $\psi$ is defined over $\Q$ and the elements of the $y_0$-fiber of $X$  
belong to a number field $k$, then, for every $i \geq 1$, the elements of the  
$A^{\circ i}(y_0)$-fiber also belong to $k$, as they are simply images of the $y_0$-fiber  
under $\psi$.

Another example illustrating Theorem \ref{t1} is given by an elliptic curve in the short Weierstrass form 
\be \l{cu} X:\, y^2=x^3+ax+b\ee defined over  $\overline\Q$. 
It is not difficult to see that, unless $a=0$,  such a curve 
does not satisfy the fourth condition of Theorem \ref{t1}, and thus does not have infinitely many $y$-fibers in a number field. On the other hand, Theorem \ref{t1} implies that curve \eqref{cu}  always has infinitely many $x$-fibers in a number field.  As a rational function $A$ in this case one can take for example the  Latt\`es map 
\be \l{llaa} A=\frac{z^4-2az^2-8bz+a^2}{4z^3+4az+4b}\ee 
associated with  the multiplication-by-$2$ endomorphism of $ X$. 

The rest of this paper is organized as follows. First, we recall a construction 
that expresses the Galois closure of the field extension $\C(X)/\C(y)$ in terms 
of fiber products. Next, we relate the fourth condition of Theorem \ref{t1} to 
the semiconjugacy relation for holomorphic maps on compact Riemann surfaces. 
Finally, we prove Theorem \ref{t1}, present several examples, and provide a 
version of Theorem \ref{t1} for finitely generated fields.

\section{Proofs}
\subsection{\l{pre} Normalizations and orbifolds}
Let $C$ be a compact  Riemann surface and 
$V:\, C\rightarrow \C\P^1$  a holomorphic map. 
A  {\it normalization} of $V$ is defined as a compact Riemann surface  $\f N_V$ together with a holomorphic Galois covering  of the 
minimal degree $\t V:\f N_V\rightarrow \C\P^1$ such that
 $\t V=V\circ H$ for some  holomorphic map $H:\,\f N_V\rightarrow C$. 
The normalization is characterized by the property that  the field extension 
$\f M(\f N_V)/\t V^*(\f M(\C\P^1))$ is isomorphic to the Galois closure 
of the extension $\f M(C)/V^*(\f M(\C\P^1))$, where $\f M(C)$ denotes the field of meromorphic functions on  a compact  Riemann surface $C$  (see e.g. \cite{des1}, Section 2.9). 
In particular, if  $X:F(x,y)=0$ is 
a curve over $\C$ and $C_X$ is its desingularization, that is,  a compact Riemann surface $C_X$ provided with holomorphic maps \linebreak $U:\, C_X\rightarrow \C\P^1$ and  $V:\, C_X\rightarrow \C\P^1$ such that the map 
\be \l{the} \theta: t\rightarrow (U(t),V(t))\subset \C\P^1\times \C\P^1\ee
is an isomorphism between $C_X$ and $X$ outside a finite set, then the condition that the Galois closure of the field extension $\C(X)/\C(y)$ has genus zero or one is equivalent to the condition that $g(\f N_V)$  equals zero or one.  

For a holomorphic map $V:\, C\rightarrow \C\P^1$ of degree  $n\geq 2$ 
its  normalization  can be described in terms of the fiber product of 
 $V$   with itself $n$ times as follows (see \cite{fried}, $\S$I.G or \cite{low}, Section 2.2).  Let $\f L^{V}$ be an algebraic variety  consisting of $n$-tuples  of $C^n$ with a common image 
under $V$, 
$$\f L^{V}= \{(x_i)\in C^n\, \vert \, V(x_1)=V(x_2)=\dots =V(x_n)\},$$
and  
 $\hat{\f L}^{V}$  a variety obtained from $\f L^{V}$ by removing components that belong to the big diagonal 
$${\Delta}^{C}:=\{(x_i)\in C^n\, \vert \, x_i=x_j \ \ {\rm for\ some} \ \ i\neq j\}$$ of  $C^n.$ 
Further, let  $\f L$ be  an arbitrary irreducible component of 
$\hat{\f L}^{V}$ and  ${\f N}\xrightarrow{\theta} \f L$ the desingularization map.  
Finally, let $\psi: {\f N}\rightarrow \C\P^1$ be a holomorphic map induced by the composition
$$  {\f N}\xrightarrow{\theta} {\f L}\xrightarrow{\pi_i}C\xrightarrow{V}{\C\P^1},$$
where  $\pi_i$ is the projection to any coordinate. In this notation, the following statement holds. 
 
\bt \l{frid} The map $\psi: {\f N}\rightarrow \C\P^1$  is the normalization of $V$. \qed
\et

 We recall that a pair $\f O=(R,\nu)$ consisting of a Riemann surface $R$ and a ramification function $\nu:R\rightarrow \mathbb N$ which takes the value $\nu(z)=1$ except at isolated points is called an {\it orbifold}. For an orbifold $\f O$ the {\it  Euler characteristic} of $\f O$ is the number
\be \l{echa} \chi(\f O)=\chi(R)+\sum_{z\in \C\P^1}\left(\frac{1}{\nu(z)}-1\right),\ee
the set of {\it singular points} of $\f O$ is the set 
$$c(\f O)=\{z_1,z_2, \dots, z_s, \dots \}=\{z\in \C\P^1 \mid \nu(z)>1\},$$ and  the {\it signature} of $\f O$ is the set 
$$\nu(\f O)=\{\nu(z_1),\nu(z_2), \dots , \nu(z_s), \dots \}.$$

If $R_1$, $R_2$ are Riemann surfaces provided with ramification functions $\nu_1,$ $\nu_2$, and 
$f:\, R_1\rightarrow R_2$ is a holomorphic branched covering map, then $f$
is called  {\it a covering map} $f:\,  \f O_1\rightarrow \f O_2$
between orbifolds
$\f O_1=(R_1,\nu_1)$ and $\f O_2=(R_2,\nu_2)$
if for any $z\in R_1$ the equality 
\be \l{us} \nu_{2}(f(z))=\nu_ {1}(z)\deg_zf\ee holds, where $\deg_zf$ stands for the local degree of $f$ at the point $z$.
%If for any $z\in R_1$ instead of equality \eqref{us} 
%the weaker condition 
%\be \l{uuss} \nu_{2}(f(z))\mid \nu_ {1}(z)\deg_zf\ee
%holds,  then the map $f$
%is called {\it a holomorphic map} $f:\,  \f O_1\rightarrow \f O_2$
%between orbifolds
%$\f O_1=(R_1,\nu_1)$ and $\f O_2=(R_2,\nu_2).$

A universal covering of an orbifold ${\f O}=(R,\nu)$
is a covering map between orbifolds  $\theta_{\f O}:\,
\tt {\f O}\rightarrow \f O$ such that $\tt R$ is simply connected and $\tt \nu(z)\equiv 1.$ 
If $\theta_{\f O}$ is such a map, then 
there exists a group $\Gamma_{\f O}$ of conformal automorphisms of $\tt R$ such that the equality 
$\theta_{\f O}(z_1)=\theta_{\f O}(z_2)$ holds for $z_1,z_2\in \tt R$ if and only if $z_1=\sigma(z_2)$ for some $\sigma\in \Gamma_{\f O}.$ A universal covering exists and 
is unique up to a conformal isomorphism of $\tt R,$
unless $\f O$ is the Riemann sphere with one ramified point, or  the Riemann sphere with two ramified points $z_1,$ $z_2$ such that $\nu(z_1)\neq \nu(z_2)$  (see \cite{fk}, 
Section IV.9.12).
Abusing  notation we will denote by $\tt {\f O}$ both the
orbifold and the  Riemann surface  $\tt R$.

With each holomorphic map $V:C\rightarrow \C\P^1$ 
one can associate in a natural way two orbifolds $\f O_1^V$
and 
$\f O_2^V$
setting $\nu_2^V(z)$  
equal to the least common multiple of local degrees of $f$ at the points 
of the preimage $V^{-1}\{z\}$, and $$\nu_1^V(z)=\nu_2^V(V(z))/\deg_zV.$$
By construction,  $V:\, \f O_1^V\rightarrow \f O_2^V$ 
is a covering map between orbifolds.  Furthermore, since a composition of covering maps  
$f:\f O_1\rightarrow \f O'$ and $g:\f O'\rightarrow \f O_2$ is a covering map   
$$g\circ f:\f O_1\rightarrow \f O_2$$ (see e.g. \cite{semi}, Corollary 4.1), 
the composition $$V\circ \theta_{\f O_1^V}: \t{\f O_1^V}\rightarrow \f O_2^V$$
is a covering map  between orbifolds, implying by the uniqueness of the universal co\-vering that 
\be \l{ravv}  \theta_{\f O_2^V}=V\circ \theta_{\f O_1^V}.\ee

For a holomorphic map $V:C\rightarrow \C\P^1$,  
the condition $g(\f  N_V)\leq 1$ can be expressed merely in terms of  $\chi({\f O}_2^V)$ as follows 
(see \cite{tame}, Lemma 2.6). 

\bl \l{ml} 
Let $C$  be a compact  Riemann surface and $V:C\rightarrow \C\P^1$ a holomorphic map. Then $g(\f  N_V)=0$ if and only if  $\chi({\f O}_2^V)> 0$, and  $g(\f  N_V)=1$ if and only if  $\chi({\f O}_2^V)= 0$. \qed 
\el

Notice that orbifolds with $\chi(\f O)\geq 0$  
and corresponding $\Gamma_{\f O}$ and $\theta_{\f O}$ can be described explicitly as follows. 
The equality $\chi(\f O)=0$ holds if and only if the signature of $\f O$ 
belongs to the list
\be \l{list}\{2,2,2,2\} \ \ \ \{3,3,3\}, \ \ \  \{2,4,4\}, \ \ \  \{2,3,6\}, \ee while $\chi(\f O)>0$  if and only if either $\f O$ is the non-ramified sphere or
 the signature of $\f O$  belongs to the  list 
 \be \l{list2} \{n,n\}, \ \ n\geq 2,  \ \ \ \{2,2,n\}, \ \ n\geq 2,  \ \ \ \{2,3,3\}, \ \ \ \{2,3,4\}, \ \ \ \{2,3,5\}.\ee

Groups $\Gamma_{\f O}\subset \Aut(\C)$ corresponding to orbifolds $\f O$ with signatures \eqref{list}  
are generated by translations of $\C$ by elements of some lattice $L\subset \C$ of rank two and the rotation $z\rightarrow  \v z,$ where $\v$ is an $n$th root of unity with $n$ equal to 2,3,4, or 6, such that  $\v L=L$.  
 Accordingly, the functions $\theta_{\f O}$ 
may be written in terms of the  corresponding
Weierstrass functions as $\wp(z),$ $\wp^{\prime }(z),$ $\wp^2(z),$  and $\wp^{\prime 2}(z)$  (see  \cite{mil2}, or \cite{fk}, 
Section IV.9.5).
Groups $\Gamma_{\f O}\subset \Aut(\C\P^1)$ corresponding to   orbifolds $\f O$ with signatures \eqref{list2} are the well-known five finite subgroups 
 $C_n,$  $D_{2n},$  $A_4,$ $S_4,$ $A_5$ of $\Aut(\C\P^1)$, and the functions $\theta_{\f O}$ are Galois coverings of $\C\P^1$ by $\C\P^1$ of degrees 
$n$, $2n,$ $12,$ $24,$  $60,$ calculated for the first time by Klein.

By Lemma \ref{ml} and equality \eqref{ravv}, for a holomorphic map $V:\, C\rightarrow \C\P^1$   the inequality $g(\f N_V)\geq 0$ holds if and only if $V$ is a ``compositional left factor'' of one of the functions  $\theta_{\f O}$ described above. Notice that this  condition is very restrictive. In particular,  %up to the change \( V \mapsto \alpha \circ V \circ \beta \), where 
%\( \alpha \) and \( \beta \) are M\"obius transformations, the list of 
rational functions \( V \) for which \( g(\f N_V) = 0 \) can be listed explicitly. 
%consists of the series 
%\[
%z^n, \ n \geq 1, \quad T_n, \ n \geq 2, \quad 
%\frac{1}{2}\left(z^n + \frac{1}{z^n}\right), \ n \geq 1,
%\]
%and a finite number of functions that can be calculated explicitly. 
On the 
other hand,  the simplest examples of rational functions with  \( g(\f N_V) = 1 \)  are 
Latt\`es maps (see \cite{gen}).

\subsection{ Semiconjugacies over $\overline\Q$}
A detailed description of the solutions to the functional equation 
\be \l{a}
A \circ V = V \circ B,
\ee
where \( A \), \( B \), and \( X \) are holomorphic maps on compact Riemann surfaces, 
was obtained in a series of papers \cite{semi}, \cite{arn}, \cite{rec}, \cite{lattes}. 
The result below partially follows from the analysis in these papers. 
However, since we require an additional conclusion concerning definability over 
\(\overline{\mathbb{Q}}\), we provide a complete proof, focusing on the existence of 
some \( A \) and \( B \) with the desired properties rather than on their full description.

\bt \l{tt} 
Let $C$  be a compact  Riemann surface and $V:C\rightarrow \C\P^1$ a holomorphic map such that 
$\chi(\f O_2^V)\geq 0$ and $c(\f O_2^V)\subset \overline{\mathbb{Q}}$. Then there exist
 holomorphic maps $A$ and $B$ of degree at least two  such that the diagram 
\be \l{i1}
\begin{CD}
C@> B>> C\\
@V {V} VV @VV {V} V\\ 
\C\P^1 @> A >> \C \P^1
\end{CD}
\ee 
commutes, the compositum of the subfields $V^*\f M(\C\P^1)$ and $B^*\f M(C)$ of $\f M(C)$ is the whole field $\f M(C)$,
 and $A$ is defined over $\overline\Q$.
 \et
\pr 
We start by observing that if the degrees of the maps $B$ and $V$ in  \eqref{i1}  are coprime, then the condition
\be \l{usl}
V^*\f M(\C\P^1) \cdot B^* \mathcal{M}(C) = \mathcal{M}(C)
\ee
is automatically satisfied. Indeed, this condition can be restated as requiring that the equalities
\[
V = \widehat{V} \circ T, \quad B = \widehat{B} \circ T,
\]
where
\[
T: C \rightarrow \t {C}, \quad \widehat{V}: \t {C} \rightarrow \C \mathbb{P}^1, 
\quad \widehat{B}: \t {C} \rightarrow C,
\]
are holomorphic maps between compact Riemann surfaces, imply that \( \deg T = 1 \). 
Since \( \deg T \) divides both \( \deg V \) and \( \deg B \), the coprimality of \( \deg V \) and \( \deg B \) ensures that this is true.

Assume first that $\chi(\f O_2^V)>0$.  In this case $\t {\f O_2^V}$ and hence $C$ is the sphere and the functions 
$\theta_{\f O_2^V}$ and $ \theta_{\f O_1^V}$ related by equality \eqref{ravv}  are rational Galois coverings.  Moreover, without loss of generality we may assume  that  $c(\f O_2^V)=\{0,\infty\}$ if  $\nu(\f O_2^V)=\{n,n\},$ $n\geq 2$, 
and 
$c(\f O_2^V)=\{0,1,\infty\}$ otherwise. Indeed,  since 
 critical values of $V$ are algebraic numbers, there exists a M\"obius map  $\nu$ with algebraic coefficients  that maps these values to $\{0,1,\infty\}$, and if $A$ satisfies the conclusions of the theorem 
 for $\nu\circ V,$ then  $\nu^{ -1} \circ A\circ \nu$ satisfies them for $V$.  
Thus, below we always will assume that $$c(\f O_2^V)\subseteq\{0,1,\infty\}.$$ 
Let us also observe that for a given rational function $V$,  we can fix  the universal covering 
$\theta_{\f O_2^V}$ in \eqref{ravv} to be any of its representatives, which are defined up to the change $\theta_{\f O_2^V}\rightarrow \theta_{\f O_2^V}\circ \mu,$ where  $\mu$ is a M\"obius transformation. Thus, below we will use such well-known representatives. 

If $\nu(\f O_2^V)=\{n,n\},$ $n\geq 2,$ then as  $\theta_{\f O_2^V}$ in \eqref{ravv} we can take $\theta_{\f O_2^V}=z^n$, and 
\eqref{ravv} implies that 
 $V=z^n\circ \mu,$ where  $\mu$ is a M\"obius transformation. Setting  now 
$$A= z^m,\ \ \ B=\mu^{-1}\circ z^m\circ \mu,$$ 
 we see that  \eqref{i1} holds and  $A$ is defined over $\overline{\Q}$.  Moreover, if $\GCD(m,n)=1,$  then  \eqref{usl} also holds. 
 
Further,  if $\nu(\f O_2^V)=\{2,2,n\}$, $n\geq 2$, then as $\theta_{\f O_2^V}$ we can take $$\theta_{\f O_2^V}=\frac{1}{2}\left(z^n+\frac{1}{z^n}\right),$$ implying by \eqref{ravv} that   either  $$V=\frac{1}{2}\left(z^n+\frac{1}{z^n}\right) \circ \mu,$$ or   $V=  T_n\circ \mu,$ 
where $T_n$ is the $n$th Chebyshev polynomial and   $\mu$ is a  M\"obius transformations (see e.g. \cite{gen}, Section 4.2). Setting now 
$$A= T_m,\ \ \ B=\mu^{-1}\circ z^m\circ \mu,$$ 
  in the first case, 
and 
$$A= T_m,\ \ \ B=\mu^{-1}\circ T_m\circ \mu$$ 
  in the second, and requiring  that $\GCD(m,n)=1,$ we see that \eqref{i1} and \eqref{usl} hold.  

Assume now that \( \nu(\f O_2^V)= \{2, 3, 4\} \). Then as $\theta_{\f O_2^V}$ in \eqref{ravv} we can take 
the function   
\be \l{th} \theta_{\f O_2^V}=-\frac{\left(z^{8}+14 z^{4}+1\right)^{3}}{108 z^{4} \left(z^{4}-1\right)^{4}},  
 \ee 
for which the corresponding group $\Gamma_{\f O}$ is generated by the M\"obius transformations  
\be \l{gru} z\rightarrow iz, \ \ \ \   z\rightarrow\frac{z+i}{z-i}.\ee
In this case, the existence of the required $A$ and $B$ can be deduced from the fact that 
 the rational function 
\be \l{alt} F=\frac{-z^5+5z}{5z^4-1},\ee 
found in \cite{genus}, satisfies the following properties: the degree of $F$ is coprime to 
\[
\deg \theta_{\f{O}_2^V} = |S_4| = 24, 
\] and 
$F$ is $\Gamma_{\f{O}_2^V}$-equivariant, 
meaning it satisfies the equality 
\be \label{homm+}
F \circ \sigma = \sigma \circ F
\ee 
for every $\sigma \in \Gamma_{\f{O}_2^V}$. 

Indeed, it follows from \eqref{homm+} that $F$ maps orbits of $\Gamma_{\f O_2^V}$ to orbits of 
$\Gamma_{\f O_2^V}$, and orbits of $\Gamma_{\f O_1^V}$ to orbits of 
$\Gamma_{\f O_1^V}$, as  $\Gamma_{\f O_1^V}$ is a subgroup of $\Gamma_{\f O_2^V}$. 
Since the orbits of $\Gamma_{\f O_2^V}$ and $\Gamma_{\f O_1^V}$ coincide 
with the fibers of $\theta_{\f O_2^V}$ and $\theta_{\f O_1^V}$, 
respectively, this implies that there exist rational functions $A$ and $B$ 
such that the diagrams 
\be \label{if} 
\begin{CD}
 \mathbb{C}\mathbb{P}^1 @> F >> \mathbb{C}\mathbb{P}^1 \\
@V \theta_{\f O_2^V} VV @VV \theta_{\f O_2^V} V \\ 
\mathbb{C}\mathbb{P}^1 @> A >> \mathbb{C}\mathbb{P}^1 
\end{CD} \ \ \ \ \
\quad\quad
\begin{CD}
 \mathbb{C}\mathbb{P}^1 @> F >> \mathbb{C}\mathbb{P}^1 \\
@V \theta_{\f O_1^V} VV @VV \theta_{\f O_1^V} V \\ 
\mathbb{C}\mathbb{P}^1 @> B >> \mathbb{C}\mathbb{P}^1
\end{CD}
\ee 
commute. Furthermore, it follows from \eqref{if} and  \eqref{ravv} that 
\[
A \circ V \circ \theta_{\f O_1^V} 
= V \circ \theta_{\f O_1^V} \circ F 
= V \circ B \circ \theta_{\f O_1^V},
\] 
implying that the diagram \eqref{i1} also 
 commutes.
In addition, since $F$ and $\theta_{\f O_2^V}$ are defined over 
$\overline{\mathbb{Q}}$, the first diagram in \eqref{if} determines 
the coefficients of $A$ through a system of linear equations over 
$\overline{\mathbb{Q}}$. Since this system has at most one solution, 
we conclude that $A$ is also defined over $\overline{\mathbb{Q}}$. 
Finally, since $$\deg V \mid \deg \theta_{\f O_2^V} = 24$$ by 
\eqref{ravv}, the degrees $\deg V$ and $\deg F$ are coprime.

If $\nu(\f{O}_2^V) = \{2, 3, 5\}$, the proof also can be deduced solely from the existence of  \(\theta_{\f{O}_2^V}\), defined over \(\overline{\mathbb{Q}}\), 
and a rational function \(F\), also defined over \(\overline{\mathbb{Q}}\), which  is $\Gamma_{\f{O}_2^V}$-equivariant and whose degree is coprime to 
\[
\deg \theta_{\f{O}_2^V} = |A_5| = 60.
\] 
Examples of such functions, as well as a description of all 
$\Gamma_{\f{O}_2^V}$-equivariant functions, can be found in 
\cite{dm}. For instance, the function $f_{11}$ of dgeree 11 given on p. 166 of \cite{dm}   is a valid 
choice for $F$.

To complete the proof in the case \( \chi(\f O_2^V) > 0 \), observe that, 
since 
$
A_4  \subset S_4,
$ 
for the remaining signature \( \nu(\f O_2^V) = \{2, 3, 3\} \), 
 we can choose the corresponding function \( \theta_{\f O_2^V} \)  to be  a compositional 
right factor of the function \eqref{th}. Consequently, the function \( F \),  defined by \eqref{alt}, remains $\Gamma_{\f O_2^V}$-equivariant for \( \nu(\f O_2^V) = \{2, 3, 3\} \), and  can  be used to complete the proof. Alternatively, since   $
A_4  \subset A_5,$  one may instead use the function $f_{11}.$

In case \( \chi(\f O_2^V) = 0 \), the proof proceeds 
in a similar way. Assume, for example, that 
\( \nu(\f O_2^V) = \{2, 2, 2, 2\} \). Then \( \theta_{\f O_2^V} \) is equal 
to the Weierstrass function \( \wp(z) \), and the group 
\( \Gamma_{\f O} \subset \operatorname{Aut}(\C) \) is generated by translations 
of \( \C \) by elements of some lattice \( L \subset \C \) of rank two, 
along with the rotation \( z \mapsto -z \). 
Since the multiplication by \( m \) on \( \C \) maps the orbits of 
\( \Gamma_{\f O_2^V} \) to the orbits of \( \Gamma_{\f O_2^V} \), and the orbits 
of \( \Gamma_{\f O_1^V} \) to the orbits of \( \Gamma_{\f O_1^V} \), arguing as 
above, we conclude that there exist holomorphic maps \( A \) and \( B \) that make 
the diagrams 
\be  
\begin{CD}
 \C @> F=mz >> \C \\
@V \theta_{\f O_2^V} VV @VV \theta_{\f O_2^V} V \\ 
\mathbb{C}\mathbb{P}^1 @> A >> \mathbb{C}\mathbb{P}^1 
\end{CD} \ \ \ \ \
\quad\quad
\begin{CD}
 \C @> F=mz >> \C \\
@V \theta_{\f O_1^V} VV @VV \theta_{\f O_1^V} V \\ 
\mathbb{C}\mathbb{P}^1 @> B >> \mathbb{C}\mathbb{P}^1
\end{CD}
\ee 
and  \eqref{i1} commutative. 
Furthermore,  since \( \deg A=\deg B = m^2 \), the degrees of \( B \) 
and \( V \) are coprime whenever \( m \) and \( \deg V \) are coprime. 

To prove that \( A \) is defined over \( \overline{\mathbb{Q}} \), we consider the algebraic 
curve \( X \) in the short Weierstrass form \eqref{cu}, parametrized by \( \wp(z) \) 
and \( \wp'(z) \). Since  \eqref{ravv} implies that the critical values of \( V \) coincide with those of \( \wp(z) \), which  are roots of the polynomial in the right part of \eqref{cu},  
it follows from the assumption $c(\f O_2^V)\subset \overline{\mathbb{Q}}$ that \( X \) is defined 
over \( \overline{\mathbb{Q}} \). Further, it is clear that multiplication by $m$ on $\mathbb{C}$ induces a 
self-rational map of $X$ of the form 
\[
\psi : (x, y) \mapsto \big(A(x), R(x, y)\big)
\] 
for some rational function $R(x, y)$. 
Finally,  since the addition operation on \( X \) is defined over the field of definition of \( X \), 
this self-rational map, and in particular the function $A$,  are also defined over \( \overline{\mathbb{Q}} \).

Further, if  
\( \nu(\f O_2^V) = \{3, 3, 3\} \), then $\theta_{\f O_2^V}=\wp^{\prime}(t)$ and the group $\Gamma_{\f O}\subset \Aut(\C)$ is  generated by translations of $\C$ by elements of a triangular lattice $L\subset \C$  and the rotation $z\rightarrow  \varepsilon z,$ where 
$\varepsilon z$ is the $3$rd root of unity. In this case, 
the proof goes similarly  with the only exception that now the required $A$ is given by the second 
coordinate of the self-rational map of $X$, induced by 
 the multiplication by \( m \) on \( \C \).  

Consider finally the case \( \nu(\f O_2^V) = \{2, 4, 4\} \); the case \( \nu(\f O_2^V) = \{2, 3, 6\} \) can be considered by an obvious modification. In this case,  \( \theta_{\f O_2^V} =\wp^2(t) \), and, as above, the multiplication by \( m \) on \( \C \) leads to holomorphic maps $B$ and $A$ that make diagram 
\eqref{i1} commutative. To show now that $A$ is defined over $\overline{\Q}$, we observe that by what is proved above $$\wp(mz)=\t A\circ \wp(z)$$ for some rational function $\t A$  defined over  $\overline{\Q}$. Therefore, it follows from 
$$\wp^2(mz)=A\circ \wp^2(z)$$ that 
$$z^2\circ \t A\circ \wp =A\circ z^2\circ \wp,$$ implying that 
$$z^2\circ \t A =A\circ z^2.$$ Since the functions $z^2$ and $\t A$ are defined over   $\overline{\Q}$, we conclude as above that the same is true for the function $A$.  
 \qed

 Our proof of Theorem \ref{tt} follows the proof of a similar existence statement 
over the field \( \mathbb{C} \) given in \cite{genus}. Notice, however, that the proof in \cite{genus} contains an erroneous claim that the function \( f_{11} \) 
commutes with all the groups \( A_4 \), \( S_4 \), and \( A_5 \). Since 
\( S_4 \not\subset A_5 \), this is not true for \( S_4 \), and a different 
\( S_4 \)-equivariant function is needed to prove the theorem for \( \nu(\f O_2^V)= \{2, 3, 4\} \). 
For example, one can use the function \eqref{alt} mentioned above, which is found in the 
same paper \cite{genus}.

\subsection{\l{dok} Proof of Theorem \ref{t1}}
To prove the implication $1)\Rightarrow 4)$, let us consider an algebraic curve $\f E$ in the space $\C^{2n}$ with coordinates 
$(x_1,y_1,x_2,y_2,\dots,x_n,y_n)$  defined by the equations 
$$F(x_1,y_1)=F(x_2,y_2)=\dots = F(x_n,y_n)=0, \ \ \ y_1=y_2=\dots =y_n.$$ 
Since the map \eqref{the} is an isomorphism off a finite set, the map 
\be \l{und} (\theta, \theta, \dots, \theta):C_X^n\rightarrow X^n\subset (\C\P^1)^{2n}\ee  
induces an isomorphism between components of the curve $\f L^{V}\subset C_X^n$ defined above and components of the curve $\f E\subset X^n$.   On the other hand,  it is easy to see that 
the first condition of the theorem implies that the image of $\hat{\f L}^{V}$ under \eqref{und} 
 has infinitely many points over $k$. Therefore, at least one of irreducible components of this image also has infinitely many points over $k$.   Thus, the fourth condition of the theorem holds by Theorem \ref{frid} and the Faltings theorem.

To prove the implication  $4)\Rightarrow 3)$, let us observe that since $X$ is defined over 
$\overline{\Q}$, critical values of $U$ and $V$ in the map \eqref{the} are algebraic. In addition, by Lemma \ref{ml}, the inequality $\chi(\f O_2^V)\geq 0$ holds. Thus,  Theorem \ref{tt} is applied to $V$ and $C=C_X$, and the  map $B$ in \eqref{i1}  obviously descends to a self-rational map $\psi$  of $X$ that makes the diagram 
\be 
\begin{CD} \l{i2} 
C_X@> B>> C_X\\
@V {(U,V)} VV @VV {(U,V)} V\\ 
X @> \psi >> X
\end{CD}
\ee 
commutative and has the form 
\be \l{by} \psi: (x,y)\rightarrow (R(x,y),A(y)),\ee where $R(x,y)$ is a  rational function over $\C$. 

Clearly,  $\psi$ maps $y$-fibers of $X$ to $y$-fibers.  Furthermore,  $R(x,y)$ can be chosen to be defined over $\overline\Q$. Indeed, if $\deg V=n$, then 
$R(x,y)$ can be represented on $X$ by some rational function $\widehat R(x,y)$ of the form 
\be \l{of} \widehat  R(x,y)=P_0(y)+P_1(y)x+P_1(y)x^2+\dots +P_{n-1}(y)x^{n-1},\ee where $P_i$, $0\leq i \leq n-1,$ are rational functions over $\C$. 
Since $X$ and $A$ are defined over $\overline\Q$,   for any algebraic point $(x_1,y_1)$ on $X$, the  points $(x_j,y_j)$, $j\geq 1$, on $ X$ defined by 
$$(x_{j+1},y_{j+1})=\left(\widehat  R(x_j,y_j),A(y_j)\right) \ \ \ j\geq 1,$$ are also algebraic. 
Thus, coefficients of  $\widehat  R(x,y)$ satisfy 
 a system of linear equations over $\overline\Q$.  Moreover, since any function of the form \eqref{of} may have only finitely many zeroes on $X$, this system has a unique solution, implying that coefficients of  $\widehat R(x,y)$ belong to $\overline\Q$. 

Finally,  it follows from \eqref{usl} 
 by the primitive element theorem that \be \l{gop2} \f M(C)=V^*\f M(\C\P^1)[h]\ee for some $h\in B^*\f M(C).$ 
As elements of  $\f M(C)$ separate 
points of $C$, 
equality \eqref{gop2} implies that for all but finitely many    
$z_0\in C$  the map $h$ takes $\deg V$ distinct values on the set $V^{-1}\{z_0\}$.  
Since $h\in B^*\f M(C),$ this yields  that for all but finitely many    
$z_0\in C$
the map $B$ takes $\deg V$ distinct values on  $V^{-1}\{z_0\}$. By \eqref{i2}, this implies in turn that 
the map induced by $\psi$ on $y$-fibers 
is generically bijective.

To prove the implication  $3)\Rightarrow 2)$, let us observe that the self-rational map $\psi:X\rightarrow X$ can be lifted to some holomorphic map $B:C_X\rightarrow C_X$ that makes the diagram 
\eqref{i2} 
commutative. Moreover, since $\psi$ sends $y$-fibers of $X$ to $y$-fibers, $B$ maps $V$-fibers to  $V$-fibers, implying that diagram \eqref{i1} commutes for some rational function $A$, and $\psi$ can be represented in  the form \eqref{by}. Furthermore, since $X$  and $\psi$ are defined over 
$\overline\Q$, there exists an infinite sequence of algebraic points $y_j$, $j\geq 1$, such that 
$A(y_j)=y_{j+1},$ which implies that $A$ is also defined over $\overline\Q$.

Given a rational map $A$, it is clear that to prove the second condition 
for a point $y_0$ it suffices to prove it for a point 
$\widehat{y}_0 = A^{\circ k}(y_0)$ for some $k \geq 1$. 
Thus, without loss of generality, we may assume that for all points $y$ 
of the forward $A$-orbit of $y_0$, the map on $y$-fibers induced by $\psi$ 
is bijective. 

Set $y_i = A^{\circ i}(y_0),$ $i \geq 0$, and let $k_0$ be a number field 
containing the coefficients of $R(x,y)$ and $A(y)$. 
By assumption, the roots 
\[
\zeta_1^i, \ \zeta_2^i, \ \dots, \ \zeta_n^i
\]
of $F(x,y_i)$, $i \geq 0$, are distinct and satisfy
\[
\zeta_1^{i+1} = R(\zeta_1^i, y_i), \quad
\zeta_2^{i+1} = R(\zeta_2^i, y_i), \quad \dots, \quad
\zeta_n^{i+1} = R(\zeta_n^i, y_i),
\]
which implies inductively that the $y_i$-fibers of $X$ lie in the field
\be \l{k} 
k = k_0(\zeta_1^0, \zeta_2^0, \dots, \zeta_n^0, y_0)
\ee
for every $i \geq 0$. 
Finally, since a finitely generated algebraic extension of a number field 
is finite, $k$ is itself a number field.

Eventually, the implication $2) \Rightarrow 1)$ is essentially obvious, 
since one can take $y_0$ to be any non-preperiodic point of $A$  
(note that, by Northcott’s theorem, such a point exists in any number field) 
and extend the field $k$ to a field containing the $A$-orbit of $y_0$. \qed

Theorem \ref{t1} implies the following corollary.  

\bc \l{coor} Let $V$ be a rational function over $\overline\Q$. Then the roots of the equation $V(x)=t$ belong to some number field $k$ for infinitely many $t\in \overline\Q$ if and only if $\chi( \f O_2^V)\geq 0.$ 
\ec
\pr To reduce Corollary \ref{coor} to Theorem \ref{t1}, it is enough to observe that if the first condition of  the corollary  is satisfied, then the corresponding values of $t$ belong to the field $\t k$ obtained from $k$ by 
 adjoining coefficients of $V$. Thus,  $\t k$ is a number field, and the curve $V(x)-y=0$ has infinitely many $y$-fibers in $\t k$. \qed

Finally, let us state an analogue of Theorem \ref{t1}, where finitely generated extensions of $\mathbb{Q}$ are considered instead of number fields.  
In accordance with the above definition, we say that a curve $X$ has infinitely many $y$-fibers in a finitely generated field if there exists such a field $k$  
and an infinite set $S \subset k$ such that, for each $y \in S$, the $y$-fiber of $X$ is contained in $k$.

\bt \l{t2}
  Let $X: F(x,y)=0$ be an  irreducible affine curve over $\C$ such that $\deg_xF\geq 2.$ Then the following conditions are equivalent. 

\begin{enumerate} [label=\arabic*)]

\item The curve $X$ has infinitely many $y$-fibers in a finitely generated field.  
\item 
There exists a rational function $A$ of degree at least two  over $\C$ satisfying the following condition:  for every  $y_0\in \C$,  
there exists a  finitely generated field $k$ such that $y$-fibers of $X$ belong to $k$ for every \( y \in \C\) that belongs to the forward orbit of $y_0$ under $A$.

\item  
There exists a self-rational map \(\psi : X \to X\) of degree at least two over $\C$,  which maps each \(y\)-fiber of \(X\) 
bijectively onto another \(y\)-fiber for all but finitely many values 
\(y \in \mathbb{C}\).

\item The Galois closure of the field extension $\C(X)/\C(y)$ has genus zero or one. 

\end{enumerate}
\et 
\pr The proof is similar to that of Theorem \ref{t1}. Specifically, the implication $1) \Rightarrow 4)$ follows from a stronger form of Faltings' theorem, where finitely generated fields are considered instead of number fields (\cite{fa}).

Furthermore, the proof of Theorem \ref{tt} clearly simplifies to yield the existence of maps $B$ and $A$ satisfying condition \eqref{i1} for any $V$ with $\chi(\f O_2^V) \geq 0$. This permits a straightforward modification of the proof of Theorem \ref{t1} to establish the implication $4) \Rightarrow 3)$.

Finally, the modifications of the proofs of the implications $2) \Rightarrow 1)$ and  $3) \Rightarrow 2)$ are straightforward since the orbit of a point $y_0$ under a rational map $A$ is contained in a finitely generated field, and the field \eqref{k} is finitely generated by construction. 
\qed

\subsection{\l{ex} Examples} 
To give an example of Corollary \ref{coor}, one can consider the orbifold $\f O$ defined by  
\be \l{db} 
\nu(0) = 2, \quad \nu(1) = 3, \quad \nu(\infty) = 4, 
\ee  
and its universal covering $\theta_{\f O}$ given by \eqref{th}. For this $\theta_{\f O}$, one can find infinitely many $t \in \Q$ such that the roots of $\theta_{\f O}(x) = t$ lie in $\Q(i)$. 
Indeed, one can check that the first diagram in \eqref{if} commutes for the functions  
$A$ and $F$ given by \eqref{aa} and \eqref{alt}.  
Thus, $F$ sends the roots of $\theta_{\f O}(z) = z_0$ to the roots of  
$\theta_{\f O}(z) = A(z_0)$ generically bijectively.  
On the other hand, if $a$ is a rational number, then all roots of  
$\theta_{\f O}(z) = \theta_{\f O}(a)$ lie in $\Q(i)$, since they form an orbit  
under the group $\Gamma_{\f O}$ generated by the Möbius transformations \eqref{gru}.

A less obvious example  is given by the rational function  \be \l{vv}  V = 3z^4 - 4z^3 \ee mentioned in the introduction, for which the orbifold $\f O$ defined by \eqref{db}  serves as $\f O_2^V$. One can check that for this function, the relation \eqref{ravv} holds for 
$$\theta_{\f O_1^V}=
\frac{\left(\frac{1}{6}(1+i){z}^{2}-\frac{i}{3}z+ \frac{1}{6}(1-i) \right)\left( {z}^{4}+2\,{z}^{3}+2\,{z}^{2}-2\,z+1 \right)}{\left( {z}^{2}+1 \right)  \left( z+1 \right)  \left( z-1 \right) z}\,, 
$$
and the diagram \eqref{i1}  commutes for  the  functions $A$ and  $B$   given by \eqref{aa} and  \eqref{bb}. 
 Thus,  $B$ sends the roots of $V(z)=z_0$ to the roots of $V(z)=A(z_0)$   generically bijectively. 
Furthermore, any algebraic curve $X : F(x, y) = 0$ of genus zero, with a  
generically one-to-one parametrization by some rational function $U$ and the  
function $V$ given by \eqref{vv}, has infinitely many $y$-fibers in a number field.  Let us observe that substituting $B(z)$ for $z$ into  
$$F(U(z), V(z)) = 0$$  
and using \eqref{i1}, one obtains  
\[
F(U \circ B(z), A \circ V(z)) = 0.
\]  
Thus, to construct for such $X$ the endomorphism $\psi$, it is enough  
to express the rational function $U \circ B$ in the form  
$U \circ B = \phi(U, V)$ for some rational function $\phi(x, y)$ over  
$\overline{\mathbb{Q}}$.

As another example, let us consider an elliptic curve given in the short Weierstrass form 
\be \l{cuu} X:y^2=x^3+ax+b, \ \ \ \ a,b\in \overline\Q,\ee where 
$4a^3+27b^2\neq 0.$ 
Clearly, the functions $U$ and $V$ in \eqref{the} have degrees two and three respectively, and,  
since $g(X)=1$, it follows from Lemma \ref{ml} that the fourth condition of Theorem \ref{t1} is satisfied 
 if and only if $c(\f O_2^V)$  belongs to the list \eqref{list}. Moreover, since $V$ has a  critical point  of multiplicity three over infinity, and the least common multiple of local degrees of a map $V$ of degree 
 three at the points of  the preimage $V^{-1}\{z\}$ cannot be equal six, the only case   $c(\f O_2^V)=\{3,3,3\}$ is possible.  The last condition is equivalent to the condition that  for some points $\pm y_0\in \C,$ the  polynomial $$x^3+ax+(b-y_0^2)$$ has a root of multiplicity two, which in turn is equivalent to the condition  $a=0.$ 
 Thus, we conclude that $X$ has infinitely many $y$-fibers in a number field
if and only if $a=0$. 
%Of course, the same result one can obtain using the parametrization of $X$ by the Weierstrass functions  by $\wp(t)$ and $\wp^{\prime }(t),$ observing that the function $\f O_1^V$ in \eqref{ravv} is non-ramified so that the condition $\chi(\f O_2^V)\geq 0$ turns out to be equivalent to a condition on the ramification of $\wp^{\prime }(t)$. 

Finally, let us observe that Theorem \ref{t1} implies that the curve \eqref{cuu} always has infinitely many $x$-fibers in a number field.  Indeed, since any extension of degree two is Galois,  for the corresponding  map $U$, we have $\widetilde U=U$ and $\f N_U=C_X$. Thus, $g(\f N_U)=g(C_X)=1$. As a function $A$ in this case we can take the Latt\`es map $A$ given by the formula \eqref{llaa}, which is defined by the commutative diagram
\be 
\begin{CD}
X@> \psi >> X\\
@V \pi VV @VV \pi  V\\ 
\C\P^1 @> A >> \C\P^1\,,
\end{CD}
\ee 
 where $\pi$ is the projection map on the $x$-coordinate and  $\psi$ is  the multiplication-by-$2$ endomorphism of $ X$.

\end{document}